# Human Behavior Algorithms for Highly Efficient Global Optimization

Da-Zheng Feng[1*], Han-Zhe Feng[1*], and Hai-Qin Zhang[2*]

[1] National Laboratory of Radar Signal Processing, Xidian University, 710071, Xi'an, P.R. China

Phone: +86 (0) 29 88201693.　Email: dzfeng@xidian.edu.cn

[2] School of Mathematics and Statistics, Xidian University, 710071, Xi'an, P.R. China

**Abstract**—The global optimization have the very extensive applications in econometrics, science and engineering. However, the global optimization for non-convex objective functions is particularly difficult since most of the existing global optimization methods depend on the local linear search algorithms that easily traps into a local point, or the random search strategies that may frequently not produce good off-springs. According to human behavior, a one-dimensional global search method in the global optimization should adopt alternating descent and ascent (up-hill and down hill) strategies. This paper proposes the human behavior algorithms (HBA) based on alternating descent and ascent approaches along a direction or multiple different directions. Very fortunately, the proposed HBA make a global optimization method have high possibility for finding a global minimum point. Several benchmark experiments test that our HBA are highly effective for solving some benchmark optimization problems.

**Index Terms**: Global optimization, Human behavior algorithm, objective function, global minimum point, local minimum point. Taylor expansion, genetic algorithm, individual, parent, offspring

## I.　INTRODUCTION

A problem in econometrics, science and engineering can be accurately or approximately formulated by a multivariate function, and its good solution is situated at the global extremum point of this function. Such function to be optimized is typically called objective function that depends on a number of variables or parameters. Because of the complexity of real-world systems, an objective function is highly nonlinear with respect to system parameters. Since the objective function naturally includes many local extremum points in which one is trapped by the conventional greedy methods such as gradient descent and Newton those [2-5, 17, 22], the challenging task for finding the global extremum may be extremely difficult.

---

* Correspondence Author: dzfeng@xidian.edu.cn; dazhengfeng@126.com
* Correspondence Author: funnyboy119@126.com
* Correspondence Author: hqzhang05@126.com

In the efforts for designing effective search means, the main research works for global optimization have been emphasis on special optimization problems and exhaustive approaches [1-4]. The various approaches can be classified as follows: (a) deterministic methods, (b) random/heuristic approaches, and (b) a mix of deterministic and random approaches. The deterministic methods include the interval methods [4], covering-based or grid-search-based methods [7]-[8], and generalized gradient methods using higher order information [9]-[13]. Conceptually, most deterministic methods adopt strategies that exploit the ingenious mechanism to find the domain of interest for the global solution [3]. Such strategies are usually improved to reduce the computational requirements by narrowing down the search regions, and discarding domains which are not likely to contain the global optima [2]-[3]. The random algorithms mainly include genetic algorithms (GA) [17], GA-typical algorithms [18], simulated annealing (SA) [19], particle swarm methods [20], ant algorithms [22], clustering methods [7], and Bayesian Methods [24]. For random methods, various random search methodologies which inherit the 'up-hill or down-hill' capability are employed to locate the global optima (see [14-15, 2]). Theoretically, probabilistic methods do also not guarantee that global optimal solutions are reached within finite search steps, even though statistically asymptotic convergence [16] can certainly be achieved until infinite iterations. The third class of methods attempts to combine merits of both deterministic and random methods to achieve efficient ways in attaining the global optimal solution.

In this paper, we propose the HBA for highly efficiently finding the global optimal solution of a nonlinear objective function. Like a human, HBA alternately use up-hill and down-hill strategies to search a global point of the nonlinear objective function. This paper is organized as follows. Section II describes HBA. In Section III, the global optimization methods are given. In section IV, the good performance of the proposed HBA is verified by several benchmark test examples. Section IV summarizes this paper. This paper uses the following notation. Boldface uppercase letters and lowercase letters denote matrices and vectors, respectively. Furthermore, $\|\bullet\|$ and $\|\bullet\|^2$ represents the Euclidean norm and its square, respectively, and $|\bullet|$ denotes the norm or magnitude operator.

## II. Human Behavior Algorithms

### A. Several Definitions and Remarks

Let a real-valued function be denoted by $f(\mathbf{x}): D \subseteq R^N \to R$, where vector $\mathbf{x} = [x_1, x_2, \cdots, x_N]^T$ is composed of $N$ variables, and $D$ is a feasible domain. Since maximization of the function $f(\mathbf{x})$ is equivalent to minimization of $-f(\mathbf{x})$, without loss of generality, we define the problem of global optimization as

$$(P): \min_{\mathbf{x}} f(\mathbf{x}), \mathbf{x} \in D. \tag{1}$$

**Definition 1:** *A point $\bar{\mathbf{x}}$ at which the gradient of $f(\mathbf{x})$ is equal to zero is called stationary point. Moreover, if there is condition $f(\bar{\mathbf{x}}) \leq f(\mathbf{x})$ or $f(\bar{\mathbf{x}}) \geq f(\mathbf{x})$, where each $\mathbf{x}$ is within infinitesimal neighborhood $\varepsilon(\bar{\mathbf{x}}) \subset D$ of $\bar{\mathbf{x}}$, then this point $\bar{\mathbf{x}}$ is referred to as the local minimum point or maximum point, respectively. Both the local minimum points and the local maximum points are uniformly called the local extremum points. Further, if Hessian matrix at point $\bar{\mathbf{x}}$ is indefinite, then the point $\bar{\mathbf{x}}$ is called saddle point.*

**Definition 2:** *A point $\bar{\mathbf{x}}$ at which $n$ partial derivatives of $f(\mathbf{x})$ with respect to $n$ variables for $1 \leq n \leq N$ are equal to zero is called $n$-order partial stationary point, where $N$ order partial stationary point represents the common stationary point. Moreover, if there is condition $f(\bar{\mathbf{x}}) \leq f(\mathbf{x})$ or $f(\bar{\mathbf{x}}) \geq f(\mathbf{x})$ only with respect to $n$ variables for $1 \leq n \leq N$, where each $\mathbf{x}$ is within infinitesimal neighborhood $\varepsilon(\bar{\mathbf{x}}) \subset D$ of $\bar{\mathbf{x}}$, then this point $\bar{\mathbf{x}}$ is referred to as the $n$ order partial local minimum point or the $n$ order partial local maximum point, respectively, where $N$ order partial local minimum point or maximum point represents the common local minimum point or the common local maximum point, respectively. An $n$ order partial local extremum point can be similarly defined.*

It is easily shown that an $n$ order partial local extremum point can not be stationary point $\bar{\mathbf{x}}$ if one of the remained $N-n$ partial derivatives of $f(\mathbf{x})$ are not equal to zeros at $\bar{\mathbf{x}}$, and may be saddle point if Hessian matrix at point $\bar{\mathbf{x}}$ is indefinite. However, an $n$ order partial local extremum point is more interesting than an arbitrary point of objective function.

**Definition 3:** *A global minimum or maximum point is defined by $\mathbf{x}^* = \arg\min_{\mathbf{x} \in D} f(\mathbf{x})$ or $\mathbf{x}^* = \arg\max_{\mathbf{x} \in D} f(\mathbf{x})$, respectively.*

We assume that problem (P) satisfies the following conditions: C1) $f(\mathbf{x})$ is continuously differentiable function; C2) the problem (P) has at least a global minimum point $\mathbf{x}^* = \arg\min_{\mathbf{x} \in D} f(\mathbf{x})$; C3) further, problem (P) includes only finite local extremum points.

It is well-known that starting from initial point $\mathbf{x}_0$, the classical optimization methods, such as gradient decent and conjugate gradient those, very easily fall in a local minimum point of $f(\mathbf{x})$. Around such local minimum point, there is an attracted domain including initial point. Starting from any initial value within this attracted domain, a classical optimization method will tend to find such local minimum point. Therefore, these classical optimization methods are especially suitable to solve various convex optimization problems, since convex optimization those have unique attracted domain. To overcome the local minimum problems, the simplest way may be to let a non-convex function be converted into the corresponding convex function by using the monotonic transformation [23] and the modern convex relax techniques. Theoretically, the monotonic transformation can convert a non-convex function into a convex function but such a convex function may include many flat regions, which can lead that many optimization methods are slow to converge. Moreover, the convex relax techniques convert only some non-convex functions into convex those, and may produce the relax gap which implies that there is significant solution difference between non-convex and convex functions. Once finding a local point and the corresponding attracted domain, some deterministic optimization methods [10-13] let the attracted domain be filed by a suitable filed function. This achieves the task for escaping the local minimum points. However, finding a suitable filed function is not an easy task.

Typically, random global optimization methods such as SA and GA escape a local minimum point and its attracted domain by adding random disturbances. However, if too small disturbances are taken, random global optimization methods do not certainly escape the attracted domain of a local minimum point. More unfortunately, too large random disturbances sometimes results in that a random global optimization method goes too far from a global minimum point. In addition, these random methods generally select arbitrary disturbances, and generate some new arbitrary points (off-springs called in GA) whose gradients are usually not equal to zero and that are also not on boundary.

***Remark 1:*** *In fact, many today global optimization methods have the search- domain-restricted problems. So is these such as SA based on statistical physical laws, and GA and ant algorithms based on biological laws. So is also those algorithms based on chemical and astronomical laws, and mathematical theorems. They find a local minimum point with a limited domain, which may be a satisfactory solution. In order to broaden the search domain, a global optimization method must take longer time. Just like our humans who take the lengthy time to recognize nature, resolving a global minimum point is long and slow. Looking around a very large domain is necessary.*

A conventional local optimization method is simply marked as LOM. An LOM belongs to one of two method classes: gradient-based class and fixed-point-based class. Its main objective is to find a local extremum point, so that it only is suitable for finding the unique extremum point of a convex function.

***Remark 2:*** *The gradient information provides a reliable guide for the conventional local optimization methods. It is well-known that water flows downwards along gradients. Similarly, a conventional local optimization method goes downwards to a local minimum point. In other words, the gradient information guides that an optimization method converges to a local extremum point. Unfortunately, for a global optimization* *problem, its gradient does not provide sufficient information that can make the state of a method tend toward a global minimum point.*

***Remark 3:*** *We note that the local maximum points and the ridges of the objective function are two of the main obstacles for finding the global minimum points. In the global optimization, it is necessary to stride over or thrill through the two obstacles. In particular, the two above obstacles can be overcome by HBA.*

**B. Human Behavior Algorithms (HBA)**

In order to efficiently search the global minimum points or the satisfied solutions, we first propose the following HBA. Let a continuously differentiable $g(x)$ with a single variable be defined within $[a,b]$. If $\frac{\partial g(a)}{\partial x} < 0$ or $\frac{\partial g(a)}{\partial x} \geq 0$, then left end point $a$ is certainly a local minimum point or a local maximum point, respectively. Hence, without losing generality, we assume that $\frac{\partial g(a)}{\partial x} < 0$. A one-dimensional local optimization method is simply marked as 1LOM. It is well-known that a good 1LOM we will adopt is the one-dimensional search method via parabolic interpolation.

The HBA consist of the following steps.

Step 1: The left end point is set as the first local maximum point.

Step 2: Starting from the left end, the next local minimum point $\tilde{x}$ is forward found by a descent 1LOM.

Step 3: Starting from the local minimum point $\tilde{x}$, the next local maximum point $\bar{x}$ is forward searched by an ascent 1LOM.

Step 4: Step 2 and Step 3 are used alternately, until the right end point is come in.

Step 5: Assume that the $\tilde{M}$ local minimum points have been found, these local minimum points are recorded as $\tilde{x}_1, \tilde{x}_2, \cdots, \tilde{x}_{\tilde{M}}$.

Since the above methods alternately use descent and ascent optimization approaches, they are intuitively called <u>a</u>lternating <u>d</u>escent and <u>a</u>scent <u>a</u>lgorithms (ADAA) or more simply referred to as HBA if not confusing. Except for finding all local minimum and maximum points along one variable, HBA tramps over hill and dale like the today fluent donkey friends or the professional mountain climbers.

***Remark 4:*** *It is worth mentioning that starting from an inter point, all the local minimum points can be obtained via the forward HBA and the back HBA. Thus, we can speak that HBA are very general algorithms.*

***Remark 5:*** *If the feasible set consists of multiple connected subsets, then HBA can be used separately to find the local minimum points within each connected subset. Finally, the (partial) local minimum points are saved and formed as a larger local minimum points list.*

Since the HBA are much complex, their simplified versions are called Taylor-HBA and followed by the following steps.

Step 1: Within each feasible connected set, $\frac{dg(x)}{dx}$ is Taylor-expanded or approximated by the polynomial with high enough order $P$ as $\frac{dg(x)}{dx} = \sum_{p=1}^{P} a_p \frac{x^p}{p!}$.

Step 2: Compute the roots of each $P$ order polynomial.

Step 3: Find all the inner points from all roots, add all the end points of the feasible connected sets, and form a large extremum point list, where an inner point implies that it belongs to feasible set; .

Step 4: Calculate all the function values of the all points with this large extremum point list, and arrange these function values in non-decreasing order.

Step 5: Store the extremum points that have the first $\tilde{M}$ smallest function values or are local minimum, and form a local minimum point list marked as $\tilde{x}_1, \tilde{x}_2, \cdots, \tilde{x}_{\tilde{M}}$ in non-decreasing order of their function values.

***Remark 6:*** *Taylor-HBA is very efficient for some special objective functions. On the basis of it, we can fast find the global minimum point of these special functions. This will be also shown via several experiments.*

## III. Global Optimization Methods

### A. GA-Based Versions

GOM based on HBA we have proposed are the modified versions of the conventional GA, and can simply call HBA if not confusing. HBA take different crossover and mutation operators, but they adopt the similar selection strategies to GA. The initial candidates (parents) with the large enough scale are randomly generated. Starting from each candidate, the corresponding local minimum point (simply called an individual like GA) with lower accuracy is sought by a LOM. Notice that aim in which the accuracy for finding a local minimum point is decreased is to make the process for performing LOM does not spend too large time. All the found individuals are formed as the initial population large enough.

**Encoding mechanism:** Fundamental to the HBA structure is the encoding mechanism for representing the optimization problem's variables. The encoding mechanism depends on the nature of the problem variables. For example, when solving for the optimal flows in a transportation problem, the variables (flows in different channels) assume continuous values, while the variables in a traveling salesman problem [21] are binary quantities representing the inclusion or exclusion of an edge in the Hamiltonian circuit. **A** large number of optimization problems have real-valued continuous variables, which is considered in this paper. A common method of encoding them uses their floating-point representation.

**Fitness function:** The objective function to be optimized provides the mechanism for evaluating each individual. However, its range of values varies from problem to problem. To maintain uniformity over various problem domains, we use the fitness function to normalize the objective function to a convenient range of 0 to 1. The normalized value of the objective function is the fitness of the individual, which the selection mechanism uses to evaluate the individuals of the population. Unfortunately, the normalized operation is usually not implemented because we do not get the range of an objective function in advance. Thus, in many references, the normalized operation is only achieved stage by stage.

**Selection:** Selection models natural survival-of-the-fittest mechanism. Fitter solutions survive, while weaker ones perish. In the HBA, a fitter individual receives a higher number of off-springs and thus has a higher chance of surviving in the subsequent generation. In the proportionate selection scheme, an individual with fitness value $f_m$ is allocated $f_m / \overline{f}$ offspring, where $\overline{f}$ is the average fitness value of the population. An individual with a fitness value higher than the average is allocated more than one offspring, while an individual with a fitness value less than the average is allocated less than one offspring. The proportionate selection scheme allocates fractional numbers of offspring to individuals. Hence the number $f_m / \overline{f}$ represents the individual's expected number of offspring. Since in the final allocation some individuals have to receive a higher number of offspring than $f_m / \overline{f}$ and some less than $f_m / \overline{f}$, allocation methods include some randomization to remove methodical allocation biases toward any particular set of offspring. The allocation technique controls the extent to which the actual allocation of offspring to individuals matches the expected number of offspring $f_m / \overline{f}$'.

***Remark 7 (good birth rule):*** *Some classical random methods produce some arbitrary off-springs. In contrast, the proposed HBA generate the off-springs that must be, at least, some partial minimum points. In fact, HBA adopt the following good birth rule: the fitness function of each offspring should be over the fitness those of its parents.*

The HBA uses the well-known roulette wheel **se**lection scheme to implement proportionate selection. Each individual is allocated a sector (slot) of a roulette wheel with the angle subtended by the sector at the center of the wheel equaling $2\pi f_m / \overline{f}$. An individual is allocated an offspring if a randomly generated number in the range 0 to $2\pi$ falls in the sector corresponding to the individual. The algorithm selects individuals in this fashion until it has generated the entire population of the next generation. Roulette wheel selection could generate large sampling errors in the sense that the final number of offspring allocated to an individual might vary significantly from the expected number. The allocated number of offspring approaches the expected number only for very large population sizes. In our selection, the individual with the highest fitness is always remained and added into the next generation.

**Crossover operators:** After selection, we are confronted with crossover, a crucial operation. Pairs of individuals are picked at random from the population to be subjected to crossover. The HBA uses the elaborately-designed crossover approach as follows. A line segment connecting pairs is given, and its feasible set is found. Along this feasible set, all the local minimum points (individuals) are searched by HBA. Such the local minimum points are usually partial. Starting from each partial local minimum point, we continuously find its local minimum point by a 1LOM, all these individuals are added into population or only their part with higher fitness are remained and participated in population.

***Remark 8:*** *If there is not any partial local minimum point between two individuals, their crossover operator is* unsuccessful, and *should be stopped. Otherwise, their crossover certainly produces a bad offspring. The above crossover always achieves good birth off-springs.*

**Mutation:** After crossover, individuals are subjected to mutation. Our mutation is also elaborately-designed. Selecting an individual from population, randomly producing a search direction, and starting from this individual along this direction, at least a local minimum point is forwards and backwards sought by HBA and added into population list.

***Remark 9:*** *Unlike the mutation used in GA, the above operator can obtain some relatively good off-springs.*

***Remark 10:*** *The above developed HBA can be seen as the GA based on good birth off-springs. The classical GA jump around any points within feasible set, which leads to some unnecessary computations. The HBA bounce only around (at least partial) local minimum point, which reduces some unnecessary computations. Hence, we expect that the HBA are more effective.*

**Competition mechanism:** The competition makes the individuals that are best suited for scanty resources survive. The HBA that search for the optimal individual can be visualized as a simultaneous competition among individuals to increase the number of their instances in the population. If we describe the optimal candidate as the individual with s high average fitness values, then the winners of the individual competitions could potentially form the optimal individual.

The genetic operators - crossover and mutation - generate, promote, and juxtapose the HBA to form optimal individuals. Crossover tends to conserve the genetic information present in the individuals to be crossed. Thus, when the individuals to be crossed are similar, its capacity to generate new candidates diminishes. Mutation is not a conservative operator and can generate radically new candidates. Selection provides the favorable bias toward candidates with higher fitness values and ensures that they increase in representation from generation to generation.

Competition is expected to makes the good parents yield good off-springs. This is not always true. Very bad off-springs can be generated when good candidates are inappropriately combined, which often appears in random crossover and mutation. The HBA select good off-springs and reject bad those as best as possible by only finding (partial) local minimum points. Off course, in order to increase the diversity of population, it may be necessary to appropriately add some bad off-springs.

*Remark 11:* *A global minimum problem may be very complex, and no algorithm can completely solve it, today. The proposed HBA can make this problem become relatively easy or its high difficulty become weak*

**D. Simplification Version Based on Alternating Coordinates Search via HBA**

It is well-known that many special problems may be very efficiently solved by some special methods. HBA for searching the partial local extreme values alternately along coordinates is just such special method as follows:

Produce $P$ initial points (parents) randomly, and then perform the following two steps until convergence.

Step 1: Passing the entire found local minimum points (parents) one by one, the partial local extreme values and the corresponding positions are sought alternately along coordinate directions by using HBA; these partial local extreme points are stored as off-springs.

Step 2: Only the $P$ (at least partial) local minimum points associated with the first $P$ highest fitness are remained and marked as $P$ parents.

*Remark 12:* *For seeking the global minimum points, it is not very necessary to select a special search direction such as steepest gradient descent direction or slowest gradient descent direction. Usually, it may be feasible to select some fixed search directions such as coordinate directions. The HBA provide a global search method along a fixed direction.*

The above algorithms are referred simply to as s-HBA. Very fortunately, although s-HBA are very simple, they work very well for finding the global point of several benchmark test examples. It is easily tested by experiments that even though only a single parent is taken, s-HBA can work efficiently.

# IV. Several Benchmark Test Examples

In this Section, from the existing references, we collect a set of benchmark objective functions. We also give the corresponding experimental results obtained by the proposed HBA.

1) Sphere model

$$f_1(\mathbf{x}) = \mathbf{x}^T \mathbf{x}, \quad \mathbf{x} = [x_1, x_2, \cdots, x_{30}]^T$$
$$-100 \leq x_n \leq 100, \quad \min f_1 = f_1(\mathbf{0}) = 0.$$

This is a variable-separable problem. Such benchmark example is too simple and is not challenging for our Taylor-HBA. It is easy to show that the Taylor-HBA can easily get the optimal solution $\mathbf{0}$ and its value 0.

2) Schwefel's problem 1

$$f_2(\mathbf{x}) = \sum_{n=1}^{30}|x_n| + \prod_{n=1}^{30}|x_n|, \quad -10 \leq x_n \leq 10$$
$$\min f_2 = f_2(\mathbf{0}).$$

This example is simpler, and is easy for our s-HBA. Readily, the s-HBA can obtain the optimal solution vector whose norm square is equal to 6.1652e-013.

3) Schwefel's problem 2

$$f_3(\mathbf{x}) = \sum_{n=1}^{30}\left(\sum_{m=1}^{n} x_m\right)^2 \quad -100 \leq x_n \leq 100$$
$$\min f_3 = f_3(\mathbf{0}) = 0.$$

This example is slightly complex, and is not difficult for our Taylor-HBA. Readily, the Taylor-HBA can get the optimal solution vector whose norm square is equal to 2.3678e-005, where the optimal value is 1.0669e-005.

4) Schwefel's problem 3

$$f_4(\mathbf{x}) = \max_n\{|x_n|, 1 \leq n \leq 30\} \quad -100 \leq x_n \leq 100$$
$$\min f_4 = f_4(\mathbf{0}).$$

This example is very simple, and is easy for our s-HBA. Readily, the s-HBA can obtain the optimal solution vector whose norm square is equal to 1.9523e-011, where the optimal value is 8.0669e-007.

5) Generalized Rosenbrock's function

$$f_5(\mathbf{x}) = \sum_{n=1}^{29}\left(100(x_{n+1} - x_n^2)^2 + (x_n - 1)^2\right) + 100(x_1 - x_{30}^2)^2 + (x_{30} - 1)^2 \quad -30 \leq x_n \leq 30$$
$$\min f_5 = f_5(1, \cdots, 1) = 0.$$

This example is slightly complex, and is not challenging for our Taylor-HBA. Readily, the Taylor-HBA can get that the optimal solution vector have elements 1.0000 and its optimal value 0.

6) Step function

$$f_6(\mathbf{x}) = \sum_{n=1}^{30}(\lfloor x_n + 0.5 \rfloor)^2 \quad -100 \leq x_n \leq 100$$
$$\min f_6 = f_6(0, \cdots, 0) = 0.$$

This example is simple, and is not challenging for our s-HBA. Readily, the s-HBA can obtain the optimal solution vector whose norm square is equal to 0, where the optimal value is 0.

7) Quartic function

$$f_7(\mathbf{x}) = \sum_{n=1}^{30} n x_n^4 \qquad -1.28 \le x_n \le 1.28$$
$$\min f_7 = f_7(0,\cdots,0) = 0.$$

This example is very simple, and is not challenging for our Taylor-HBA. Readily, the Taylor-HBA can get the optimal solution vector whose norm square is equal to 0, where the optimal value is 0.

8) Generalized Schwefel's problem 3

$$f_8(\mathbf{x}) = \sum_{n=1}^{30}\left(x_m \sin\left(\sqrt{|x_m|}\right)\right) = \sum_{m=1}^{30} g_8(x_m)$$
$$-500 \le x_n \le 500$$
$$\min f_8 = f_8(420.9687,\cdots,420.9687) = -12569.5.$$

This example is complex and has too many local minimum points. Fig.1 shows the curves of $g_8(x)$, where all local minimum and maximum points are denoted by symbol "*". Fortunately, since it is a variable-separable problem and may be not very challenging for our s-HBA. Readily, the s-HBA can obtain that the optimal solution vector have elements -420.9687, where the optimal value is -1.2569e+004.

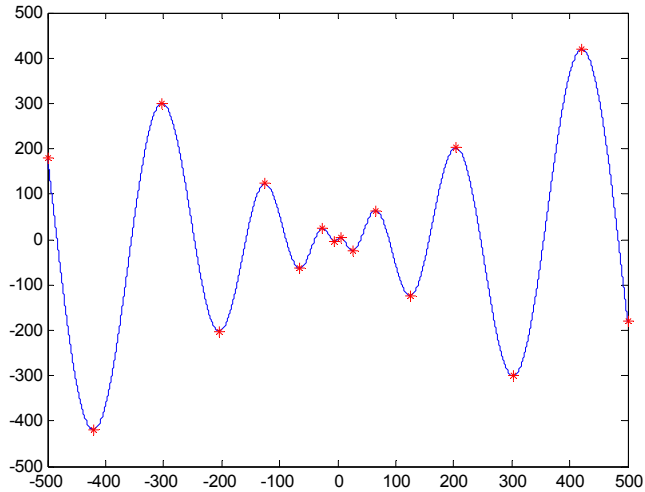

Fig. 1. Curves of $g_8(x)$ is shown, where all local minimum and maximum points are gotten by s-HBA, and denoted as "*".

9) Generalized Rastrigin's function

$$f_9(\mathbf{x}) = \sum_{n=1}^{30}\left[x_n^2 - 10\cos(2\pi x_n) + 10\right] = \sum_{n=1}^{30} g_9(x_n)$$
$$-5.12 \le x_n \le 5.12$$
$$\min f_9 = f_9(0,\cdots,0) = 0.$$

This example is complex and has too many local minimum points. Fig.2 shows the curves of $g_9(x)$, where all local minimum and maximum points are denoted by symbol "*". Fortunately, since it is a variable-separable problem and may be not very challenging for our s-HBA. Readily, the s-HBA can get that the optimal solution vector have elements 0.5797e-013, where the optimal value is 0.

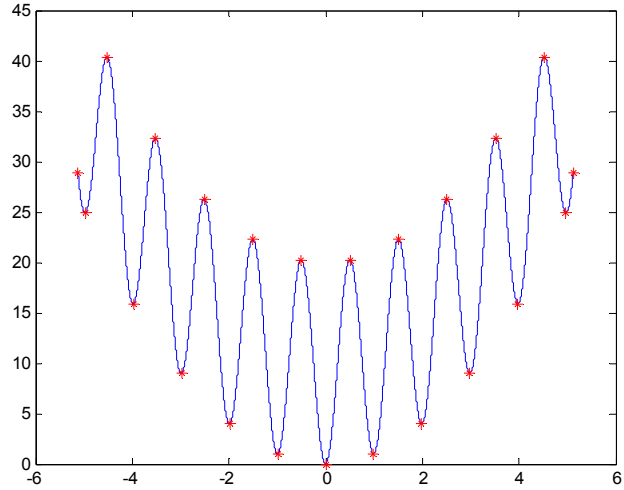

Fig. 2. Curves of $g_9(x)$ is shown, where its all local minimum and maximum points are gotten by s-HBA, and denoted as "*".

10) Ackley's function

$$f_{10}(\mathbf{x}) = -20\exp\left(-0.2\sqrt{\frac{1}{30}\sum_{n=1}^{30}x_n^2}\right) - \exp\left(\frac{1}{30}\sum_{n=1}^{30}\cos(2\pi x_n)\right) + 20 + e$$

$$-32 \leq x_n \leq 32, \quad \min f_{10} = f_{10}(0,\cdots,0) = 0.$$

This example seems to be complex and has too many local minimum points. Fig.3 shows a slice curve of $f_{10}(\mathbf{x})$. Fortunately, since the attraction domain of each local point is very small, it may be not very challenging for our s-HBA. Readily, the s-HBA can obtain that the optimal solution vector whose norm square equals to 2.0440e-023, where the optimal value is 1.1546e-014.

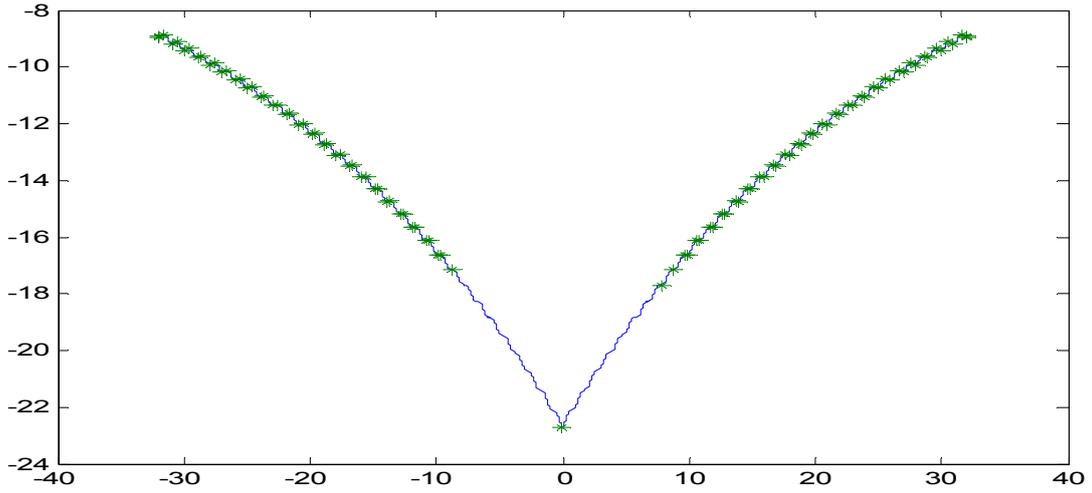

Fig.3. A slice curve of $f_{10}(\mathbf{x})$ is shown, where its all local minimum and maximum points are gotten by s-HBA, and denoted as "*".

11) Generalized Griewank's function

$$f_{11}(\mathbf{x}) = \frac{1}{4000}\sum_{n=1}^{30} x_n^2 - \prod_{n=1}^{N} \cos\left(\frac{x_n}{\sqrt{n}}\right) + 1$$
$$-600 \le x_n \le 600, \quad \min f_{11} = f_{11}(0,\cdots,0) = 0.$$

This example seems to be complex and has too many local minimum points. Fortunately, since the attraction domain of each local point is very small, it may be not very challenging for our s-HBA. Readily, the s-HBA can get that the optimal solution vector whose norm square equals to 2.6014e-014, where the optimal value is -1.1728e-006. Fig.4 shows a slice curve of $f_{11}(\mathbf{x})$ near the optimal point.

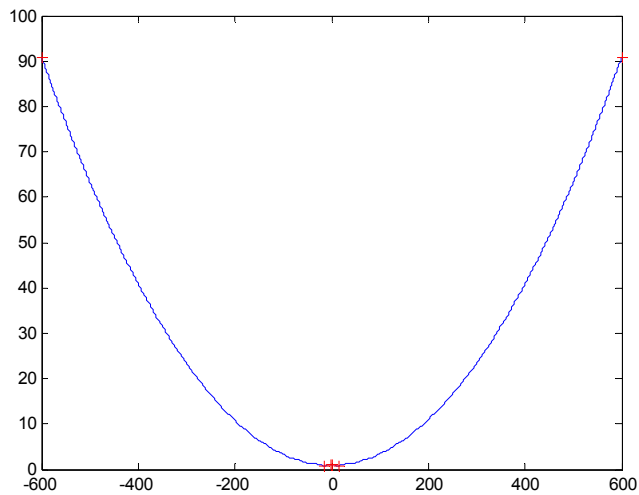

Fig.4. A slice curve of $f_{11}(\mathbf{x})$ is shown near the optimal point, where its all local minimum and maximum

12) First Modified Griewank's function

$$f_{12}(\mathbf{x}) = \frac{1}{4000}\sum_{n=1}^{30} x_n^2 - \prod_{n=1}^{N}\left[\frac{2}{3} + \frac{1}{3}\cos\left(\frac{x_n}{\sqrt{n}}\right)\right] + 1$$

$$-100 \leq x_n \leq 100, \quad \min f_{12} = f_{12}(0,\cdots,0) = 0.$$

This example seems to be complex and has too many local minimum points. Fig.5 shows a slice curve of $f_{12}(\mathbf{x})$. Fortunately, it may be not very challenging for our s-HBA. Readily, the s-HBA can obtain that the optimal solution vector whose norm square equals to 2.6014e-014, where the optimal value is -1.1728e-006.

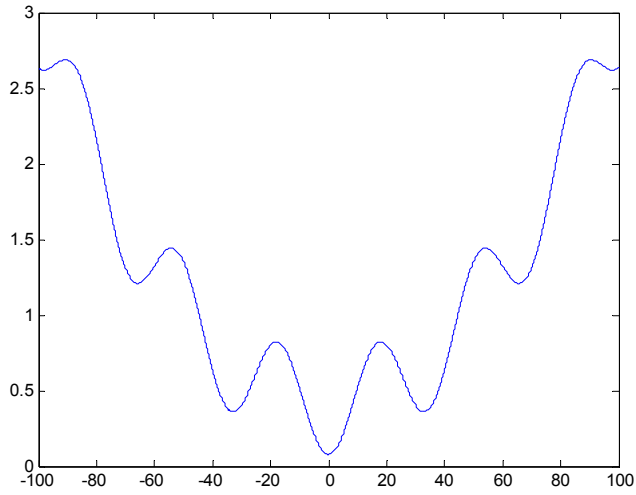

Fig.5. A slice curve of $f_{12}(\mathbf{x})$ near the optimal point, where its all local minimum and maximum points are not denoted for clarity.

13) Second modified Griewank's function

$$f_{13}(\mathbf{x}) = \frac{1}{4000}\sum_{n=1}^{30} x_n^2 - \prod_{n=1}^{N}\ln\left[2 + \cos\left(\frac{x_n}{\sqrt{n}}\right)\right] + (\ln 3)^{30}$$

$$-600 \leq x_n \leq 600, \quad \min f_{13} = f_{13}(0,\cdots,0) = 0.$$

This example is complex and has too many local minimum points. Fig.6 shows a slice curve of $f_{13}(\mathbf{x}) - (\ln 3)^{30}$. Fortunately, although it is difficult for many methods, the s-HBA can get that the optimal solution vector whose norm square equals to 3.0156e-012, where the optimal value is 5.6133e-013.

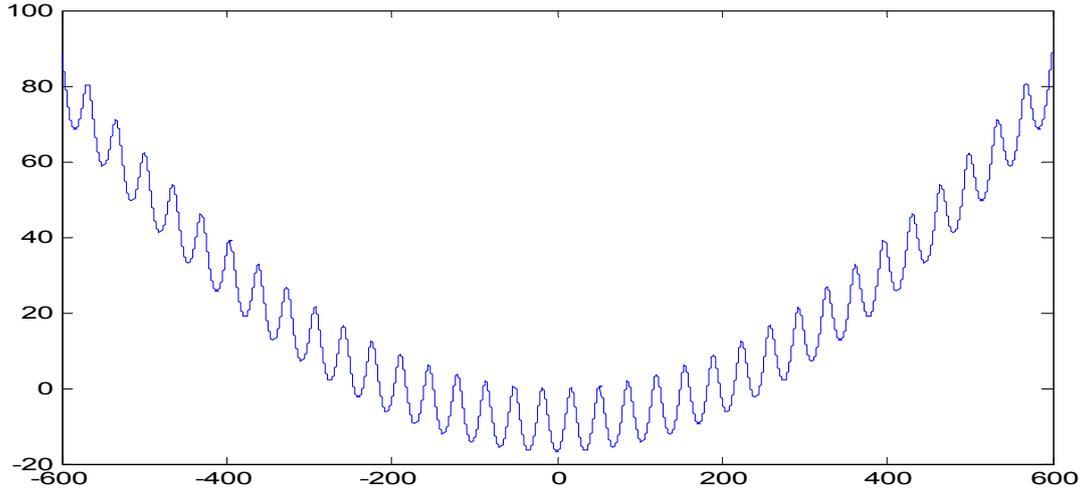

Fig.6. A slice curve of $f_{13}(\mathbf{x}) - (\ln 3)^{30}$ near the optimal point, where its all local minimum and maximum points are not denoted for clarity.

14) Third modified Griewank's function

$$f_{14}(\mathbf{x}) = \frac{1}{4000}\sum_{n=1}^{30}\sum_{m=1}^{30} x_n a_{n,m} x_m - \prod_{n=1}^{30} \ln\left[2+\cos(x_n)\right] + (\ln 3)^{30}$$
$$-150 \leq x_n \leq 150, \quad \min f_{14} = f_{14}(0,\cdots,0) = 0$$
where $a_{n,m} = 1$ for $m \neq n$, and $a_{n,n} = n$ for $m = n$.

This example seems to be complex and has too many local minimum points. It is extremely difficult for many methods. However, our s-HBA can deal with it. The s-HBA can obtain the optimal solution vector whose the norm square is 9.5211e-015, where the optimal value is -1.7764e-014. Fig.7 shows a slice curve of $f_{14}(\mathbf{x})$ near the optimal solution.

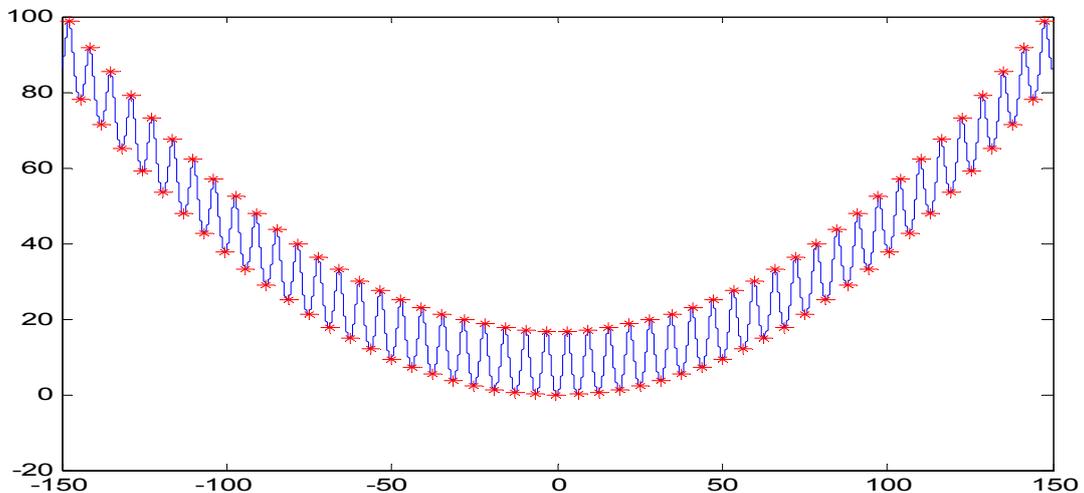

Fig.7. A slice curve of $f_{14}(\mathbf{x}) - (\ln 3)^{30}$ near the optimal point is shown, where its all local minimum and maximum points are gotten by s-HBA, and denoted as "*".

15) Fourth modified Griewank's function

$$f_{15}(\mathbf{x}) = \frac{1}{10000} \sum_{n=1}^{30} x_n^2 - \sum_{n=1}^{30} \ln\left[2 + \cos\left(\frac{x_n}{\sqrt{n}}\right)\right] + 30\ln 3$$

$$= \sum_{n=1}^{30} g_{15,n}(x_n)$$

$$-200 \leq x_n \leq 200, \quad \min f_{15} = f_{15}(0,\cdots,0) = 0.$$

This example seems to be complex and has too many local minimum points. Fig.8 shows the curve of $g_{15,16}(x)$. It is very difficult for many methods. Fortunately, since it is a separable problem, our s-HBA can deal relatively easily with it. The s-HBA can get that the optimal solution vector whose the norm square is 4.6593e-012, where the optimal value is 4.3299e-014.

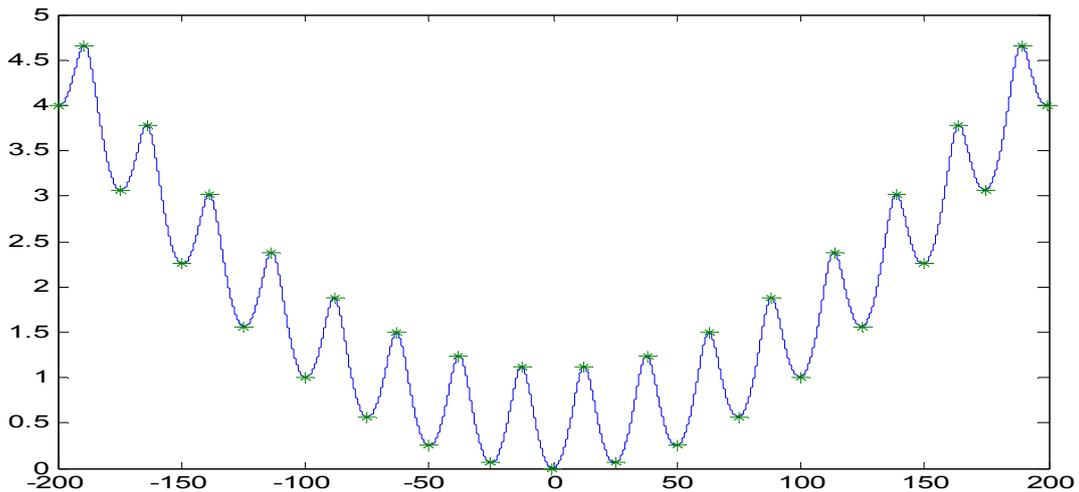

Fig.8. The curve of $g_{15,16}(x)$ is shown near the optimal point, where its all local minimum and maximum points are gotten by s-HBA, denoted as "+".

16) Fifth modified Griewank's function

$$f_{16}(\mathbf{x}) = \frac{1}{40000} \sum_{n=1}^{30} \sum_{m=1}^{30} x_n a_{n,m} x_m - \sum_{n=1}^{30} \ln\left[2 + \cos(x_n)\right] + 30\ln 3$$

$$-150 \leq x_n \leq 150, \quad \min f_{16} = f_{16}(0,\cdots,0) = 0$$

where $a_{n,m} = 1$ for $m \neq n$, and $a_{n,n} = n$ for $m = n$.

This example seems to be complex and has too many local minimum points. It is extremely difficult for many methods. However, our s-HBA can deal with it. The s-HBA can obtain the optimal solution vector whose the norm square is 2.3686e-021, where the optimal value is -1.4211e-014. Fig.9 shows a slice curve of $f_{16}(\mathbf{x})$ near the optimal solution.

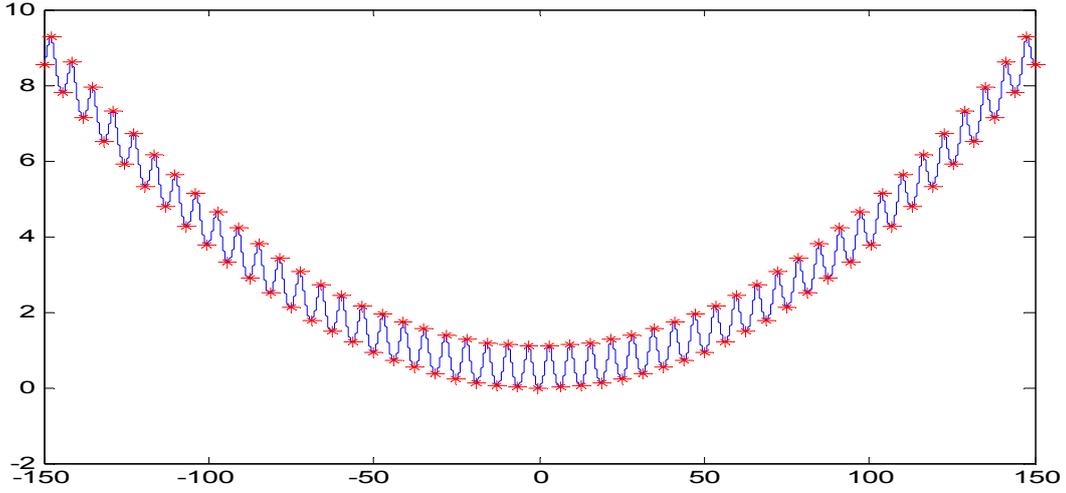

Fig.9. A slice curve of $f_{16}(\mathbf{x})$ near the optimal point is shown, where its all local minimum and maximum points are gotten by s-HBA, and denoted as "*".

17) Generalized Penalized function 1

$$f_{17}(\mathbf{x}) = \frac{\pi}{30}\left\{10\sin^2(\pi y_1) + \sum_{n=1}^{29}(y_n-1)^2[1+10\sin^2(\pi y_{n+1})] + (y_{30}-1)^2\right\}$$

$$= \sum_{n=1}^{30} u(x_n, 10, 100, 4)$$

$$-50 \le x_n \le 50, \quad \min f_{17} = f_{17}(-1,\cdots,-1) = 0$$

$$\text{where } u(x_n, a, k, m) = \begin{cases} k(x_n-a)^m, & x_n > a \\ 0, & -a \le x_n \le a \\ k(-x_n-a)^m & x_n < -a \end{cases}$$

$$y_n = 1 + (x_n+1)/4$$

This example is simpler, and is not challenging for our s-HBA. Readily, the s-HBA can get the optimal solution vector whose entries are equal to -1.0000, where the optimal value is 1.5705e-032.

18) Generalized Penalized function 2

$$f_{18}(\mathbf{x}) = 0.1\left\{\sin^2(3\pi x_1) + \sum_{n=1}^{29}(x_n-1)^2[1+10\sin^2(3\pi x_{n+1})] + (x_{30}-1)^2[1+10\sin^2(3\pi x_{30})]\right\}$$

$$= \sum_{n=1}^{30} u(x_n, 5, 100, 4)$$

$$-50 \le x_n \le 50, \quad \min f_{18} = f_{18}(1,\cdots,1) = 0$$

This example is simpler, and is not challenging for our s-HBA. Readily, the s-HBA can obtain the optimal solution vector whose entries are equal to 1.0000, where the optimal value is 3.8114e-020.

19) Shekel's function

$$f_{19}(x_1, x_2) = \left[\frac{1}{500} + \sum_{m=1}^{25} \frac{1}{m + \sum_{n=1}^{2}(x_n - a_{n,m})^6}\right]^{-1}$$

$$-66 \leq x_n \leq 66, \quad \min f_{19} = f_{19}(-32, -32) \approx 1.$$

where $(a_{n,m}) = \begin{bmatrix} -32 & -16 & 0 & 16 & 32 & -32 & \cdots & 0 & 16 & 32 \\ -32 & -16 & 0 & 16 & 32 & -32 & \cdots & 0 & 16 & 32 \end{bmatrix}$

This example is simpler, and is not challenging for our s-HBA. Readily, the s-HBA can get the optimal solution vector whose entries are equal to -32.0450, where a better optimal value is found to be 0.7302. Fig.10 shows a slice curve of $f_{19}(\mathbf{x})$ near the optimal value.

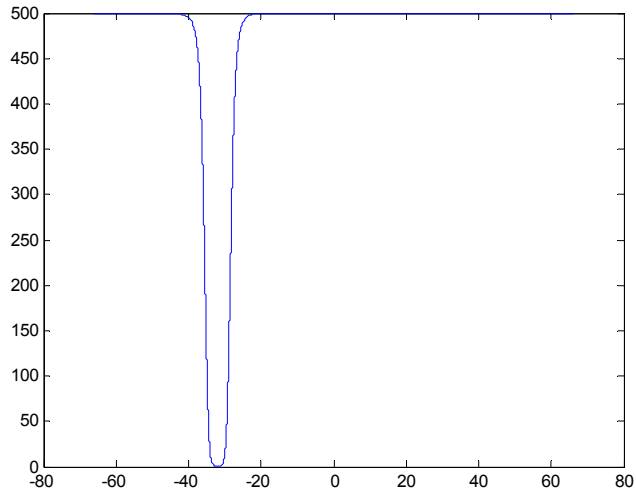

Fig.10. A slice curve of $f_{19}(\mathbf{x})$ near the optimal point.

20) Six-Hump Camel-Back function

$$f_{20}(x_1, x_2) = 4x_1^2 - 2.1x_1^4 + \frac{1}{3}x_1^6 + x_1 x_2 - 4x_2^2 + 4x_2^4$$

$$-5 \leq x_n \leq 5$$

$$\mathbf{x}_{\min} = [0.08983, -0.7126]^T, [-0.08983, 0.7126]^T$$

$$\min f_{20} = -1.0316285.$$

This example is simpler, and is not challenging for our HBA. Readily, the HBA can obtain the optimal solution vector whose two entries are equal to -0.0898 and 0.7127, where the optimal value is -1.0316. Fig.11 shows a slice curve of $f_{20}(\mathbf{x})$ near the optimal value.

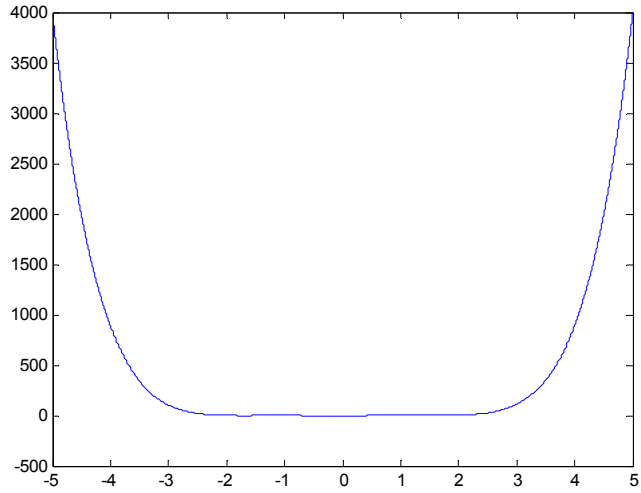

Fig.11. A slice curve of $f_{20}(\mathbf{x})$ near the optimal point.

21) Branin function

$$f_{21}(x_1, x_2) = \left(x_2 - \frac{5.1}{4\pi^2}x_1^2 + \frac{5}{\pi}x_1 - 6\right)^2 + 10\left(1 - \frac{1}{8\pi}\right)\cos x_1 + 10$$
$$-5 \le x_1 \le 10, \quad 0 \le x_2 \le 15$$
$$\mathbf{x}_{\min} = [-3.142, 12.275]^T, [9.425, 2.425]^T$$
$$\min f_{21} = 0.398.$$

This example is simpler, and is not challenging for our s-HBA. Readily, the s-HBA can get the optimal solution vector whose two entries are equal to 3.1416 and 5.0055, which is a new optimal point, where the better optimal value is 0.3979. Fig.12 shows a slice curve of $f_{21}(\mathbf{x})$ near the optimal value.

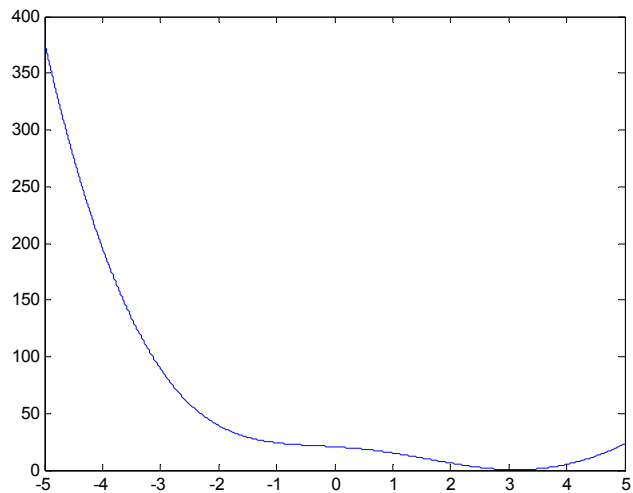

Fig.12. A slice curve of $f_{21}(\mathbf{x})$ near the optimal point.

22) Goldstein-Price function

$$f_{22}(x_1,x_2) = \left[1 + (x_1 + x_2 + 1)^2(19 - 14x_1 + 3x_1^2 + 6x_1x_2 - 14x_2 + 3x_2^2)\right]$$
$$\times \left[30 + (2x_1 - 3x_2)^2(18 - 32x_1 + 12x_1^2 - 36x_1x_2 + 48x_2 + 27x_2^2)\right].$$
$$-2 \leq x_n \leq 2, \quad \min f_{22} = f_{22}(0,-1) = 3$$

This example is simpler, and is not challenging for our HBA. Readily, the HBA can obtain the same optimal solution vector whose two entries are equal to -0.0000 and 1.0000, which is a new optimal point, where the optimal value is 3.0000.

23) Hartman's Function family

Let

$$\mathbf{A} = \begin{bmatrix} 10.0 & 3.00 & 17.0 & 3.50 & 1.70 & 8.00 \\ 0.05 & 10.0 & 17.0 & 0.10 & 8.00 & 14.0 \\ 3.00 & 3.50 & 1.70 & 10.0 & 17.0 & 8.00 \\ 17.0 & 8.00 & 0.05 & 10.0 & 0.10 & 14.0 \end{bmatrix}$$

$$\mathbf{c} = \begin{bmatrix} 1.0 & 1.2 & 3 & 3.2 \end{bmatrix}^T$$

$$\mathbf{P} = \begin{bmatrix} 0.1312 & 0.1696 & 0.5569 & 0.0124 & 0.8283 & 0.5886 \\ 0.2329 & 0.4135 & 0.8307 & 0.3736 & 0.1004 & 0.9991 \\ 0.2348 & 0.1415 & 0.3522 & 0.2883 & 0.3047 & 0.6650 \\ 0.4047 & 0.8828 & 0.8732 & 0.5743 & 0.1091 & 0.0381 \end{bmatrix}$$

$$f(\mathbf{x},n) = -\sum_{i=1}^{4} c_i \exp\left[\sum_{j=1}^{n} a_{i,j}(x_j - p_{i,j})^2\right]$$

Then

$$f_{23}(\mathbf{x}) = f(\mathbf{x},3)$$
$$f_{24}(\mathbf{x}) = f(\mathbf{x},6)$$
$$0 \leq x_n \leq 1,$$
$$f_{23}(0.114, 0.556, 0.852) = -3.86$$
$$f_{24}(0.201, 0.150, 0.477, 0.275, 0.311, 0.657) = -3.32$$

The two examples are simpler, and are not challenging for our s-HBA. For the first example, the s-HBA can get the new better optimal solution vector whose three entries are equal to 0.1858, 0.2001 and 0.5594, where the better optimal value is -3.9552. For the second example, the s-HBA can obtain the similar optimal solution vector whose six entries are equal to 0.2017, 0.1468, 0.4767, 0.2753, 0.3117 and 0.6573, where the optimal value is -3.3220.

# V. Conclusions

We propose HBA to imitate human's up-hill and down-hill behaviors. The HBA can be applied to significantly improve many global optimization methods. A set of benchmark objective functions are simulated, and the obtained experimental results show that HBA are very effective.